\documentclass[12pt]{article}
\usepackage{amsmath}
\usepackage{amssymb}
\usepackage{amsbsy}
\usepackage{fancyhdr}
\usepackage{xcolor}
\usepackage{amsfonts}
\newtheorem{lem}{Lemma}[section]%
\newtheorem{theorem}[lem]{Theorem}%
\newtheorem{prop}[lem]{Proposition}%
\newtheorem{claim}{Claim}
\newtheorem{case}{Case}

\def\a{\alpha} \def\b{\beta} \def\g{\gamma} \def\d{\delta} 
 \def\s{\sigma}

\def\G{\Gamma}

 \def\lg{\langle} \def\rg{\rangle}
\def\nd{\mathrel{\bigm|\kern-.7em/}}

\def\f{\noindent}

\def\z{\mathbb{Z}}
\def\PSL{\hbox{\rm PSL}}

\def\GL{\hbox{\rm GL}}
\def\PGL{\hbox{\rm PGL}}

\def\Aut{\hbox{\rm Aut\,}}

\def\Ker{\hbox{\rm Ker }}

\def\Cay{\hbox{\rm Cay }}

\def\mod{\hbox{\rm mod }}

\def\F{\hbox{\rm F}}
\def\C{\hbox{\rm C}}

\def\BiCay{{\rm BiCay}}

\def\R{{\mathcal R}}
\def\L{{\mathcal L}}
\def\H{{\mathcal H}}

\def\demo{{\bf Proof}\hskip10pt}

\def\mz{{\mathbb Z}}
\def\qed{\hskip10pt $\Box$\vspace{3mm}}

\begin{document}
\title{Cubic edge-transitive bi-Cayley graphs over inner-abelian $p$-groups
\author
{Yan-Li Qin, Jin-Xin Zhou\\
{\small Mathematics, Beijing Jiaotong University, Beijing 100044, P.R. China}\\
{\small{Email}:\ 15118413@bjtu.edu.cn, jxzhou@bjtu.edu.cn}
}}

\date{}
\maketitle

\begin{abstract}
A graph is said to be a {\em bi-Cayley graph} over a group $H$ if it admits $H$ as a group of automorphisms
acting semiregularly on its vertices with two orbits. A non-abelian group is
called an {\em inner-abelian group} if all of its proper subgroups are abelian.
In this paper, we complete the classification of connected cubic edge-transitive
bi-Cayley graphs over inner-abelian $p$-groups for an odd prime $p$.

\bigskip

\end{abstract}

\section{Introduction}

Throughout this paper, we denote by $\z_n$ the cyclic group of order $n$ and
by $\z_n^*$ the multiplicative group of $\z_n$ consisting of numbers coprime to $n$.
All groups are assumed to be finite, and all graphs are assumed to be finite, connected, simple and undirected.
Let $\G=(V(\G), E(\G))$ be a graph with vertex set $V(\G)$, and edge set $E(\G)$. Denote by $\Aut(\G)$ the full automorphism group of $\G$. For $u, v\in V(\G)$, denote by $\{u, v\}$ the edge incident to $u$ and $v$ in $\G$. For a graph $\G$, if $\Aut(\G)$ is transitive on $V(\G)$ or $E(\G)$, then $\G$ is said to be {\em vertex-transitive} or {\em edge-transitive}, respectively. An arc-transitive graph is also called a symmetric graph. 

Let $G$ be a permutation group on a set $\Omega$ and take $\a\in \Omega$.
The stabilizer $G_{\a}$ of $\a$ in $G$ is the subgroup of $G$ fixing the point $\a$.
The group $G$ is said to be {\em semiregualr} on $\Omega$ if $G_{\a}=1$ for every $\a\in \Omega$ and
{\em regular} if $G$ is transitive and semiregular.

A graph is said to be a {\em bi-Cayley graph} over a group $H$ if
it admits $H$ as a semiregular automorphism group with two orbits (Bi-Cayley graph is sometimes called {\em semi-Cayley graph}).
Note that every bi-Cayley graph admits the following concrete realization. Given a group $H$, let $\R$, $\L$ and $S$ be subsets of $H$ such that $\R^{-1}=\R$, $\L^{-1}=\L$ and $\R\cup \L$ does not contain the identity element of $H$. The {\em bi-Cayley graph} over $H$ relative to the triple $(\R, \L, S)$, denoted by BiCay($H, \R, \L, S$), is the graph having vertex set the union $H_{0}\cup H_{1}$ of two copies of $H$, and edges of the form  $\{h_{0},~(xh)_{0}\}$, $\{h_{1},~(yh)_{1}\}$ and $\{h_{0},~(zh)_{1}\}$ with $x\in \R, y\in \L, z\in S$ and $h_0\in H_0$, $h_1\in H_1$ representing a given $h\in H$.
Let $\G=\BiCay(H, \R, \L, S)$. For $g\in H$, define a permutation $R(g)$ on the vertices of $\G$ by the rule
$$h_i^{R(g)}=(hg)_i, \forall i\in\mz_2, h\in H.$$
Then $R(H)=\{R(g)\ |\ g\in H\}$ is a semiregular subgroup of $\Aut(\G)$ which is isomorphic to $H$ and has $H_0$ and $H_1$ as its two orbits.
When $R(H)$ is normal in $\Aut(\G)$, the bi-Cayley graph $\G=\BiCay(H,\R,\L,S)$ is called a {\em normal bi-Cayley graph} over $H$ (see \cite{Zhouaut}).
A bi-Cayley graph $\G=\BiCay(H,\R,\L,S)$ is called {\em normal edge-transitive} if $N_{\Aut(\G)}(R(H))$ is transitive on the edge-set of $\G$ (see \cite{Zhouaut}).

There are many important graphs which can be constructed as bi-Cayley graphs. For example,
the Petersen graph is a bi-Cayley graph over a cyclic group of order $5$. Another interesting
bi-Cayley graph is the Gray graph~\cite{Bouwer1} which is a bi-Cayley graph over a metacyclic $3$-group of order $27$.
One more example of bi-Cayley graph is the Hoffman-Singleton graph~\cite{HS-graph} which is a bi-Cayley graph over
an elementary abelian group of order $25$. We note that all of these graphs are bi-Cayley graphs over a $p$-group.
Inspired by this, we are naturally led to investigate the bi-Cayley graphs over a $p$-group.

In \cite{Zhouaut}, a characterization is given of cubic edge-transitive bi-Cayley graphs over a $2$-group.
A next natural step would be studying cubic edge-transitive bi-Cayley graphs over a $p$-group, where $p$ is an odd prime.
Due to Zhou et al.'s work in \cite{Zhoucubic} about the classification of cubic vertex-transitive abelian bi-Cayley graphs, we
may assume the $p$-group in question is non-abelian. As the beginning of this program, in \cite{QYL} we prove that every cubic edge-transitive bi-Cayley graph over a $p$-group is normal whenever $p>7$, and
moreover, it is shown that a cubic edge-transitive bi-Cayley graph over a metacyclic $p$-group exists only when $p=3$,
and cubic edge-transitive bi-Cayley graphs over a metacyclic $p$-group are normal except the Gray graph. Recall that a non-abelian group is
called an {\em inner-abelian group} if all of its proper subgroups are abelian. We note that the Gray graph \cite{Bouwer1},
the smallest cubic semisymmetric graph, is isomorphic to $\BiCay(H,\emptyset,\emptyset,\{1, a, a^2b\})$, where $H$ is
the following inner-abelian metacyclic group of order $27$
\[\lg a,b\ |\ a^9=b^3=1, b^{-1}ab=a^4\rg.\] In \cite{QYL}, a complete classification
is given of cubic edge-transitive bi-Cayley graphs over an inner-abelian metacyclic $p$-group.

In this paper, we shall complete the classification of cubic edge-transitive bi-Cayley graphs over any inner-abelian $p$-group. By \cite{Redei} or \cite[Lemma~65.2]{p-group2}, for every odd prime $p$, an inner-abelian non-metacyclic $p$-group is isomorphic to
the following group:
\begin{equation}\label{inner}
\H_{p, t, s}=\lg a, b, c\ |\ a^{p^t}=b^{p^s}=c^p=1, [a, b]=c, [c, a]=[c, b]=1\rg (t\geq s\geq 1).
\end{equation}

Now we define a family of cubic bi-Cayley graphs over $\H_{p,t,s}$. If $t=s$, then take $k=0$, while if $t>s$, take $k\in\mz_{p^{t-s}}^*$ such that $k^2-k+1\equiv 0\ (\mod p^{t-s})$. Let
\begin{equation}\label{graph}
\Sigma_{p, t, s, k}=\BiCay(\H_{p, t, s},\emptyset, \emptyset, \{1, a, ba^k\}).
\end{equation}
It will be shown in Lemma~\ref{Isomorphic} that for any two distinct admissible integers $k_1, k_2$, the graphs $\Sigma_{p, t, s, k_1}$ and $\Sigma_{p, t, s, k_2}$ are isomorphic. So the graph $\Sigma_{p, t, s, k}$ is independent of the choice of $k$, and we denote by $\Sigma_{p, t, s}$
the graph $\Sigma_{p, t, s, k}$.

Before stating our main result, we introduce some symmetry properties of graphs. An {\em $s$-arc}, $s\geq 1$, in a graph $\G$ is an ordered $(s+1)$-tuple $(v_0,v_1, \ldots ,v_{s-1},v_s)$ of vertices of $\G$ such that $v_{i-1}$ is adjacent to $v_i$ for $1\leq i\leq s$, and $v_{i-1} \neq v_{i+1}$ for $1 \leq i \leq s-1$, and a $1$-arc is usually called an {\em arc}. A graph $\G$ is said to be {\em $s$-arc-transitive} if $\Aut(\G)$ is
transitive on the set of $s$-arcs in $\G$. An $s$-arc-transitive graph is said to be {\em $s$-transitive} if it is not $(s+1)$-arc-transitive. In particular, $0$-arc-transitive means {\em vertex-transitive}, and $1$-arc-transitive means {\em arc-transitive} or {\em symmetric}. A subgroup $G$ of $\Aut(\G)$ is {\em $s$-arc-regular} if for any two $s$-arcs of $\G$, there is a unique element $g\in G$ mapping one to the other, and $\G$ is said to be {\em $s$-arc-regular} if $\Aut(\G)$ is $s$-arc-regular.
It is well known that, in the cubic case, an $s$-transitive graph is $s$-arc-regular.

\begin{theorem}\label{classify-all}
Let $\G$ be a connected cubic edge-transitive bi-Cayley graph over $\H_{p,s,t}$.
Then $\G\cong \Sigma_{p, t, s}$. Furthermore, the following hold:
\begin{enumerate}
  \item [{\rm (1)}]\ $\Sigma_{3, 2, 1}$ is $3$-arc-regular;
  \item [{\rm (2)}]\ $\Sigma_{p,t,s}$ is $2$-arc-regular if $t=s$;
  \item [{\rm (3)}]\ $\Sigma_{3,t,s}$ is $2$-arc-regular if $t=s+1$, and $(t, s)\neq(2, 1)$;
  \item [{\rm (4)}]\ $\Sigma_{p,t,s}$ is $1$-arc-regular if $p^{t-s}>3$.
\end{enumerate}
\end{theorem}




We shall close this section by introducing some notation which will be used in this paper.
For a finite group $G$, the full automorphism group, the center, the derived subgroup and the Frattini subgroup of $G$ will be denoted by $\Aut(G)$, $Z(G)$, $G'$ and $\Phi(G)$, respectively. For $x, y\in G$, denote by $o(x)$ the order of $x$ and by $[x, y]$ the commutator $x^{-1}y^{-1}xy$. For a subgroup $H$ of $G$,
denote by $C_G(H)$ the centralizer of $H$ in $G$ and by $N_G(H)$ the
normalizer of $H$ in $G$. For two groups $M$ and $N$, $N\rtimes M$ denotes a semidirect product of $N$ by $M$.

\section{Some basic properties of the group $\H_{p,t,s}$}

In this section, we will give some properties of the group $\H_{p,t,s}$ (given in Equation~(\ref{inner})).

\begin{lem}\label{comput2}
Let $H=\H_{p,t,s}$. Then the following hold:
\begin{enumerate}
  \item [{\rm (1)}]\ For any $i\in \z_p^t$, we have $a^ib=ba^ic^i$.

  \item [{\rm (2)}]\ $H'=\lg c\rg\cong\z_p$.

  \item [{\rm (3)}]\ For any $x,y\in H$, we have $(xy)^p=x^py^p$.

  \item [{\rm (4)}]\ For any $x,y\in H$, if $o(x)=o(a)=p^t, o(y)=o(b)=p^s$ and $H=\lg x, y\rg$,
  then $H$ has an automorphism taking $(a, b)$ to $(x, y)$.

  \item [{\rm (5)}]\ Every maximal subgroup of $H$ is one of the following groups:
  \[\begin{array}{l}

 \lg ab^j, b^p, c\rg=\lg ab^j\rg\times\lg b^p\rg\times\lg c\rg\cong\z_{p^{t}}\times\z_{p^{s-1}}\times\z_p(j\in\mz_p), \\

 \lg a^p,b, c\rg=\lg a^p\rg\times\lg b\rg\times\lg c\rg\cong\z_{p^{t-1}}\times\z_{p^{s}}\times\z_p.\\
\end{array}
\]
\end{enumerate}
\end{lem}

\f\demo For (1), for any $i\in\mz_{p^t}$, since $[a, b]=c$ and $[c, a]=1$, we have $b^{-1}ab=ac$ and $ac=ca$, and then $b^{-1}a^ib=(b^{-1}ab)^i=(ac)^i=a^ic^i$.
It follows that $a^ib=ba^ic^i$, and so (1) holds.

From \cite[Lemma~65.2]{p-group2}, we have the items (2) and (3).

For (4), assume that $H=\lg x, y\rg$, and $o(x)=o(a), o(y)=o(b)$. Let $z=[x, y]$. Then $z\neq 1$, and then by (2),
we have $H'=\lg z\rg=\lg c\rg$. It follows that $z^p=1$ and $[z,x]=[z,y]=1$. Consequently, $x$ and $y$ have the
same relations as do $a$ and $b$. Therefore, $H$ has an automorphism taking $(a, b)$ to $(x, y)$.

For (5), let $M$ be a maximal subgroup of $H$. As $H$ is a $2$-generator group, we have $H/\Phi(H)=\lg a\Phi(H)\rg\times\lg b\Phi(H)\rg\cong\z_p\times\z_p$. Clearly, $M/\Phi(H)$ is a subgroup of $H/\Phi(H)$ of order $p$, so $M/\Phi(H)=\lg ab^j\Phi(H)\rg$
or $\lg b\Phi(H)\rg$ for some $j\in\mz_p$.
Note that $\Phi(H)=\lg a^p, b^p, c\rg$ is contained in the center of $H$. It follows that $M$ is one of the following groups:
  \begin{align*}
 & \lg ab^j, b^p, c\rg=\lg ab^j\rg\times\lg b^p\rg\times\lg c\rg\cong\z_{p^{t}}\times\z_{p^{s-1}}\times\z_p(j\in\mz_p), \\
 & \lg a^p,b, c\rg=\lg a^p\rg\times\lg b\rg\times\lg c\rg\cong\z_{p^{t-1}}\times\z_{p^{s}}\times\z_p.\\
\end{align*}
This proves (5). \hfill\qed


\section{The isomorphisms of $\Sigma_{p,t,s,k}$}
The goal of this section is to prove the graph $\Sigma_{p,t,s,k}$ is independent on the choice of $k$.
By the definition, if $t=s$, then $k=0$, and so for any given group $\H_{p,t,s}$, we only have one graph.
So we only need to consider the case when $t>s$. We first restate an easily proved result about bi-Cayley graphs.

\begin{prop}{\rm~\cite[Lemma~3.1]{Zhoucubic}}\label{bicayley}
Let $\G=\BiCay(H,\R,\L,S)$ be a connected bi-Cayley graph over a group $H$. Then the following hold:
\begin{enumerate}
\item[$(1)$] $H$ is generated by $\R\cup \L\cup S$.
\item[$(2)$] Up to graph isomorphism, $S$ can be chosen to contain the identity of $H$.
\item[$(3)$] For any automorphism $\alpha$ of $H$, $\BiCay(H, \R, \L, S)\cong \BiCay(H, \R^{\alpha}, \L^{\alpha}, S^{\alpha})$.
\end{enumerate}
\end{prop}

\begin{lem}\label{Isomorphic}
Suppose that $t>s$ and $k_1, k_2\in\mz_{p^{t-s}}^*$ are two distinct solutions of the equation $k^2-k+1\equiv0\ (\mod p^{t-s})$. Then
$\Sigma_{p, t, s, k_1} \cong \Sigma_{p, t, s, k_2}$.
\end{lem}

\f\demo Recall that
$$\H_{p, t, s}=\lg a, b, c\ |\ a^{p^t}=b^{p^s}=c^p=1, [a, b]=c, [c, a]=[c, b]=1\rg, $$
and $$\Sigma_{p, t, s, k_i}=\BiCay(\H_{p, t, s},\emptyset, \emptyset, T_{i}), {\rm where}~ T_{i}=\{1, a, ba^{k_i}\}\ {\rm with} \ i=1,2. \ \ \  $$

We first show that there exists an automorphism $\b$ of $\H_{p,t,s}$ which sends $(a,b)$ to $(ba^{k_2}, a(ba^{k_2})^{-k_1})$.
It is easy to see that $ba^{k_2}, a(ba^{k_2})^{-k_1}$ generate $\H_{p,t,s}$. By Lemma~\ref{comput2}~(4), it suffices to show
that $o(a)=o(ba^{k_2})$ and $o(b)=o(a(ba^{k_2})^{-k_1})$. Since $k_2\in\mz_{p^{t-s}}^*$, from Lemma~\ref{comput2}~(3) it follows that
$o(a)=o(ba^{k_2})$. By Lemma~\ref{comput2}~(3), we have $(a(ba^{k_2})^{-k_1})^{p^s}=a^{p^s}(b^{p^s}a^{k_2p^s})^{-k_1}=(a^{p^s})^{1-k_1k_2}$.
Since $k_1, k_2\in\mz_{p^{t-s}}^*$ satisfy $k^2-k+1\equiv0\ (\mod p^{t-s})$, it follows that $-k_1, -k_2$ are two elements of $\mz_{p^{t-s}}^*$ of order $3$. Since $\mz_{p^{t-s}}^*$ is cyclic, we have $k_1k_2\equiv 1\ (\mod p^{t-s})$. Consequently, $(a^{p^s})^{1-k_1k_2}=1$ and so $o(a(ba^{k_2})^{-k_1})=o(b)$.


Now we know that $\H_{p, t, s}$ has an automorphism $\b$ taking $(a, b)$ to $(ba^{k_2}, a(ba^{k_2})^{-k_1})$.
Moreover,
$$T_{k_1}^{\b}=\{1, a, ba^{k_1}\}^{\b}=\{1, ba^{k_2}, a(ba^{k_2})^{-k_1}\cdot (ba^{k_2})^{k_1}\}=\{1, ba^{k_2}, a\}=T_{k_2}. $$
By Proposition~\ref{bicayley}~(3),
we have \[\Sigma_{p, t, s, k_1}=\BiCay(\H_{p, t, s},\emptyset, \emptyset, T_{1})\cong\BiCay(\H_{p, t, s},\emptyset, \emptyset, T_{2})=\Sigma_{p, t, s, k_2},\] as required.
\hfill\qed

\section{The automorphisms of $\Sigma_{p,t,s}$}

The topic of this section is the automorphisms of $\Sigma_{p,t,s}$.

\subsection{Preliminaries}

In this subsection, we give some preliminary results.
Let $\G$ be a connected graph with an edge-transitive group $G$ of automorphisms
and let $N$ be a normal subgroup of $G$. The {\em quotient graph} $\G_N$ of $\G$ relative to $N$ is defined as the graph
with vertices the orbits of $N$ on $V(\G)$ and with two orbits
adjacent if there exists an edge in $\G$ between the vertices lying in
those two orbits. Below we introduce two propositions, of which the
first is a special case of \cite[Theorem~9]{VTgraph}.

\begin{prop}\label{3orbits}
Let $\G$ be a cubic graph and let $G\leq \Aut(\G)$ be arc-transitive on $\G$.
Then $G$ is an $s$-arc-regular subgroup of $\Aut(\G)$ for some integer $s$.
If $N\unlhd G$ has more than two orbits in $V(\G)$, then $N$ is semiregular on $V(\G)$,
$\G_N$ is a cubic symmetric graph with $G/N$ as an $s$-arc-regular subgroup of automorphisms.
\end{prop}

The next proposition is a special case of \cite[Lemma~3.2]{6p2}.

\begin{prop}\label{intransitive}
Let $\G$ be a cubic graph and let $G\leq \Aut(\G)$ be transitive on $E(\G)$ but intransitive on $V(\G)$.
Then $\G$ is a bipartite graph with two partition sets, say $V_0$ and $V_1$.
If $N \trianglelefteq G$ is intransitive on each of $V_0$ and $V_1$, then $N$ is semiregular on $V(\G)$,
$\G_N$ is a cubic graph with $G/N$ as an edge- but not vertex-transitive group of automorphisms.
\end{prop}

The following result gives an upper bound of the order of the vertex-stabilizer of cubic edge-transitive graphs.

\begin{prop}{\rm~\cite[Proposition~8]{2p3}}\label{Av}
Let $\G$ be a connected cubic edge-transitive graph and let $G\leq\Aut(\G)$ be transitive on the edges of $\G$. For any $v\in V(\G)$, the stabilizer $G_v$ has order $2^r\cdot 3$ with $r\geq 0$.
\end{prop}

The next three propositions are about cubic edge-transitive bi-Cayley graphs over a $p$-group.

\begin{prop}{\rm~\cite[Lemma~4.1]{QYL}}\label{N}
Let $\G$ be a connected cubic edge-transitive graph of order $2p^n$ with $p$ an odd prime and $n\geq 2$.
Let $G\leq\Aut(\G)$ be transitive on the edges of $\G$.
Then any minimal normal subgroup of $G$ is an elementary abelian $p$-group.
\end{prop}

\begin{prop}{\rm~\cite[Lemma~4.2]{QYL}}\label{11}
Let $p\geq 5$ be a prime and let $\G$ be a connected cubic edge-transitive graph of order $2p^n$ with $n\geq 1$.
Let $A=\Aut(\G)$ and let $H$ be a Sylow $p$-subgroup of $A$. Then $\G$ is a bi-Cayley graph over $H$, and moreover,
if $p\geq 11$, then $\G$ is a normal bi-Cayley graph over $H$.
\end{prop}

\begin{prop}{\rm~\cite[Lemma~4.3]{QYL}}\label{Q}
Let $\G$ be a connected cubic edge-transitive graph of order $2p^n$ with $p=5$ or $7$ and $n\geq 2$.
Let $Q=O_p(A)$ be the maximal normal $p$-subgroup of $A=\Aut(\G)$. Then $|Q|=p^n$ or $p^{n-1}$.
\end{prop}

\subsection{Normality of cubic edge-transitive bi-Cayley graphs over $\H_{p,t,s}$}

The following lemma determines the normality of cubic edge-transitive bi-Cayley graphs over $\H_{p,t,s}$.

\begin{lem}\label{57}
Let $\G$ be a connected cubic edge-transitive bi-Cayley graph over $\H_{p,t,s}$.
If $p=3$, then $\G$ is normal edge-transitive. If $p>3$, then $\G$ is normal.
\end{lem}

\f\demo Let $A=\Aut(\G)$ and let $P$ be a Sylow $p$-subgroup of $A$ such that $R(H)\leq P$. Let $H=\H_{p,t,s}$, and let $|H|=p^n$ with $n={t+s+1}$.
If $p=3$, then by Proposition~\ref{Av}, we have $|A|=3^{n+1}\cdot 2^r$ with $r\geq 0$. This implies that $|P|=3|R(H)|$, and so $|P_{1_0}|=|P_{1_1}|=3$. Thus, $P$ is transitive on the edges of $\G$. Clearly, $R(H)\unlhd P$. This implies that $\G$ is normal edge-transitive.

Suppose now $p>3$. Then $R(H)$ is a Sylow $p$-subgroup of $A$. Suppose to the contrary that $R(H)$ is not normal in $A$. By Proposition~\ref{11}, we have $p=5$ or $7$.
Let $N$ be the maximal normal $p$-subgroup of $A$. Then $N\leq R(H)$, and by Proposition~\ref{Q}, we have $|R(H): N|=p$. By Propositions~\ref{3orbits} and \ref{intransitive}, the quotient graph $\G_N$ is a cubic graph of order $2p$ with $A/N$ as an edge-transitive automorphism group. By \cite{Sym768, Semisym768}, if $p=5$, then $\G_N$ is the Petersen graph, and if $p=7$, then $\G_N$ is the Heawood graph. Since $A/N$ is transitive on the edges of $\G_N$ and $R(H)/N$ is non-normal in $A/N$, it follows that
$$
\begin{array}{ll}
A_5\lesssim A/N\lesssim S_5, & {\rm if }\ p=5;\\
\PSL(2,7)\lesssim A/N\lesssim \PGL(2,7), & {\rm if }\ p=7.
\end{array}
$$
Let $B/N$ be the socle of $A/N$. Then $B/N$ is also edge-transitive on $\G_N$, and so $B$ is also edge-transitive on $\G$.
Let $C=C_B(N)$. Then $C/(C\cap N)\cong CN/N \unlhd B/N$. Since $B/N$ is non-abelian simple, one has $CN/N=1$ or $B/N$.

Suppose first that $CN/N=1$. Then $C\leq N$, and so $C=C\cap N=C_N(N)=Z(N)$.
Since $R(H)$ is inner-abelian, we have $N$ is abelian, and so $C=Z(N)=N$.
Recall that $|R(H): N|=p$. Then $N$ is a maximal subgroup of $R(H)$.
By Lemma~\ref{comput2}~(5), we have $N\cong \z_{p^t}\times\z_{p^{s-1}}\times\z_p$ or $\z_{p^{t-1}}\times\z_{p^s}\times\z_p$.
Let $\mho_1(N)=\{x^p\ |\ x\in N\}$ and $M=(R(H))'\mho_1(N)$.
Then $\mho_1(N)\cong \z_{p^{t-1}}\times\z_{p^{s-2}}$ or $\z_{p^{t-2}}\times\z_{p^{s-1}}$. Moreover, $M$ is characteristic in $N$ and $N/M\cong \mathbb{Z}_p\times\mathbb{Z}_p$. It implies that each element $g$ of $B$ induces an automorphism of $N/M$, denote by $\s(g)$. Consider the map $\varphi: B \rightarrow\Aut(N/M)$ with $\varphi (g)=\s(g)$ for any $g\in B$. It is easy to check that $\varphi$ is a homomorphism. Letting $\Ker\varphi$ be the kernel of $\varphi$, we have $\Ker\varphi=C=N$. It follows that $B/N\lesssim\Aut(N/M)\cong\GL(2, p)$. This forces that either $A_5\leq \GL(2,5)$ with $p=5$, or $\PSL(2,7)\leq \GL(2,7)$ with $p=7$. However, each of these can not happen by Magma{~\cite{Magma}}, a contradiction.

Suppose now that $CN/N=B/N$. Since $C\cap N=Z(N)$, we have $1<C\cap N\leq Z(C)$.
Clearly, $Z(C)/(C\cap N) \unlhd C/(C\cap N)\cong CN/N$. Since $CN/N=B/N$ is non-abelian simple, $Z(C)/C\cap N$ must be trivial.
Thus $C\cap N=Z(C)$, and hence $B/N=CN/N\cong C/C\cap N=C/Z(C)$. If $C=C'$, then $Z(C)$ is a subgroup of the Schur multiplier of $B/N$. However, the Schur multiplier of $A_5$ or $\PSL(2,7)$ is $\mz_2$, a contradiction. Thus, $C\neq C'$.
Since $C/Z(C)$ is non-abelian simple, one has $C/Z(C)=(C/Z(C))'=C'Z(C)/Z(C)\cong C'/(C'\cap Z(C))$, and then we have $C=C'Z(C)$. It follows that $C''=C'$. Clearly, $C'\cap Z(C)\leq Z(C')$, and $Z(C')/(C'\cap Z(C)) \unlhd C'/(C'\cap Z(C))$. Since $C'/(C'\cap Z(C))\cong C/Z(C)$ and since $C/Z(C)$ is non-abelian simple, it follows that $Z(C')/(C'\cap Z(C))$ is trivial, and so $Z(C')=C'\cap Z(C)$.
As $C/(C\cap N)\cong CN/N$ is non-abelian, we have $C/(C\cap N)=(C/(C\cap N))'=(C/Z(C))'\cong C'/(C'\cap Z(C))=C'/Z(C')$.
Since $C'=C''$, $Z(C')$ is a subgroup of the Schur multiplier of $CN/N$. However, the Schur multiplier of $A_5$ or $\PSL(2,7)$ is $\mz_2$, forcing that $Z(C')\cong\mz_2$. This is impossible because $Z(C')=C'\cap Z(C)\leq C\cap N$ is a $p$-subgroup. Thus $R(H)\unlhd A$, as required. \medskip
\hfill\qed

\subsection{Automorphisms of $\Sigma_{p,t,s}$}
We first collect several results about the automorphisms of the bi-Cayley graph $\G={\rm BiCay}(H, \R, \L, S)$.
Recall that for each $g\in H$, $R(g)$ is a permutation on $V(\G)$ defined by the rule
\begin{equation}\label{1}
h_{i}^{R(g)}=(hg)_{i},~~~\forall i\in \mathbb{Z}_{2},~h,~g\in H,
\end{equation}
and $R(H)=\{R(g)\ |\ g\in H\}\leq\Aut(\G)$.
For an automorphism $\a$ of $H$ and $x,y,g\in H$, define two permutations on $V(\G)=H_0\cup H_1$ as following:
\begin{equation}\label{2}
\begin{array}{ll}
\d_{\a,x,y}:& h_0\mapsto (xh^{\a})_1, ~h_1\mapsto (yh^{\a})_0, ~\forall h\in H,\\
\s_{\a,g}:& h_0\mapsto (h^\a)_0, ~h_1\mapsto (gh^{\a})_1, ~\forall h\in H.\\
\end{array}
\end{equation}
Set \begin{equation}\label{3}
\begin{array}{lll}
{\rm I}&=& \{\d_{\a,x,y}\ |\ \a\in\Aut(H)\ s.t.\ \R^\a=x^{-1}\L x, ~\L^\a=y^{-1}\R y, ~S^\a=y^{-1}S^{-1}x\},\\
{\rm F} &=&\{ \s_{\a,g}\ |\ \a\in\Aut(H)\ s.t.\ \R^\a=\R, ~\L^\a=g^{-1}\L g, ~S^\a=g^{-1}S\}.
\end{array}
\end{equation}

\begin{prop}{\rm~\cite[Theorem~3.4]{Zhouaut}}\label{bicayleyaut}
Let $\Gamma=\BiCay(H, \R, \L, S)$ be a connected bi-Cayley graph over the group $H$.
Then $N_{\Aut(\Gamma)}(R(H))=R(H)\rtimes F$ if $I=\emptyset$ and
$N_{\Aut(\Gamma)}(R(H))=R(H)\langle F, \delta_{\a,x,y}\rangle$ if
$I\neq\emptyset$ and $\delta_{\a,x,y}\in I$. Furthermore, for any $\delta_{\a,x,y}\in I$, we have the following:
\begin{enumerate}
\item[$(1)$] $\langle R(H), \delta_{\a,x,y}\rangle$ acts transitively on $V(\Gamma)$;
\item[$(2)$] if $\a$ has order $2$ and $x=y=1$, then $\Gamma$ is isomorphic to the Cayley graph $\Cay(\bar{H},~\R\cup \alpha S)$, where $\bar{H}=H\rtimes \langle\a\rangle$.
\end{enumerate}
\end{prop}

\begin{lem}\label{Example2}
The graph $\Sigma_{p, t, s}$ is symmetric.
\end{lem}

\f\demo Recall that
$$\H_{p, t, s}=\lg a, b, c\ |\ a^{p^t}=b^{p^s}=c^p=1, [a, b]=c, [c, a]=[c, b]=1\rg, $$
and $$\Sigma_{p, t, s}=\BiCay(\H_{p, t, s},\emptyset, \emptyset, \{1, a, ba^{k}\}), $$
where if $t=s$, then $k=0$, and if $t>s$, then $k\in\mz_{p^{t-s}}^*$ satisfies $k^2-k+1\equiv0\ (\mod p^{t-s})$.

We first prove the following two claims.\medskip

\f{\bf Claim 1} $\H_{p, t, s}$ has an automorphism $\a$ mapping $a, b$ to $a^{-1}ba^k, a^{-1}(a^{-1}ba^k)^{-k}$, respectively.

By definition, if $t=s$ then $k=0$, and by Lemma~\ref{comput2}~(4), we can obtain Claim~1. Let $t>s$. Then $k^2-k+1\equiv0\ (\mod p^{t-s})$. Let $x=a^{-1}ba^k$ and $y=a^{-1}(a^{-1}ba^k)^{-k}$. Note that $(yx^k)^{-1}=a$ and $(yx^k)^{-1}x(yx^k)^{k}=b$. This implies that $\lg x, y\rg=\lg a, b\rg=\H_{p, t, s}$.

By Lemma~\ref{comput2} (1), we have $x=a^{-1}ba^k=ba^{k-1}c^{-1}$. Since $k^2-k+1\equiv0\ (\mod p^{t-s})$, we have $(k-1, p)=1$.
By Lemma~\ref{comput2} (3), we have $o(x)=o(a)=p^t$.
Since $p^{t-s}\mid k^2-k+1$, again by Lemma~\ref{comput2} (3),
$$y^{p^s}=(a^{-1}(a^{-1}ba^k)^{-k})^{p^s}=a^{-p^s}(a^{-p^s}b^{p^s}a^{kp^s})^{-k}=(a^{-p^s})^{k^2-k+1}=1, $$
and so $o(y)=o(b)=p^s$.
By Lemma~\ref{comput2} (4), $\H_{p, t, s}$ has an automorphism taking $(a, b)$ to $(x, y)$, as claimed.\medskip

\f{\bf Claim 2.} $\H_{p, t, s}$ has an automorphism $\b$ mapping $a, b$ to $a^{-1}, a^{-k}b^{-1}a^k$, respectively.

Let $u=a^{-1}$ and $v=a^{-k}b^{-1}a^k$. Clearly, $\lg u,v\rg=\lg a,b\rg=\H_{p, t, s}$ and $o(u)=p^t$.
Note that
$$v^{p^s}=(a^{-k}b^{-1}a^k)^{p^s}=a^{-kp^s}b^{-p^s}a^{kp^s}=1. $$
So $o(v)=o(b)=p^s$.
By Lemma~\ref{comput2} (4), $\H_{p, t, s}$ has an automorphism taking $(a, b)$ to $(u, v)$, as claimed.\medskip

Now we are ready to finish the proof of our lemma. Set $T=\{1, a, ba^k\}$. By Claim~1, there exists $\a\in\Aut(\H_{p, t, s})$ such that $a^{\a}=a^{-1}ba^k$ and $b^{\a}=a^{-1}(a^{-1}ba^k)^{-k}$. Then
$$a^{-1}T=a^{-1}\{1, a, ba^k\}=\{a^{-1}, 1, a^{-1}ba^k\},$$
$$T^\a=\{1, a, ba^k\}^{\a}=\{1, a^{-1}ba^k, a^{-1}(a^{-1}ba^k)^{-k}\cdot (a^{-1}ba^k)^k\}=\{1, a^{-1}ba^k, a^{-1}\}. $$
Thus $T^{\a}=a^{-1}T$. By Proposition~\ref{bicayleyaut}, $\s_{\a, a}$ is an automorphism of $\Sigma_{p, t, s}$ fixing $1_0$ and cyclically permutating the three neighbors of $1_0$.
Set $B=R(\H_{p, t, s})\rtimes\lg \s_{\a,a}\rg$. Then $B$ acts transitively on the edges of $\Sigma_{p, t, s}$.

By Claim~2, there exists $\b\in\Aut(\H_{p, t, s})$ such that $a^{\b}=a^{-1}$ and $b^{\b}=a^{-k}b^{-1}a^k$.
Then $$T^{\b}=\{1, a, ba^k\}^{\b}=\{1, a^{-1}, a^{-k}b^{-1}a^k\cdot a^{-k}\}=\{1, a^{-1}, a^{-k}b^{-1}\}=T^{-1}.$$
By Proposition~\ref{bicayleyaut}, $\d_{\b,1,1}$ is an automorphism of $\Sigma_{p, t, s}$ swapping $1_0$ and $1_1$.
Thus, $\Sigma_{p, t, s}$ is vertex-transitive, and so $\Sigma_{p, t, s}$ is symmetric.
\hfill\qed

\begin{theorem}\label{s-regular}
One of the following holds.
\begin{enumerate}
  \item [{\rm (1)}]\ $\Sigma_{3,2,1}$ is $3$-arc-regular;
  \item [{\rm (2)}]\ $\Sigma_{p,t,s}$ is $2$-arc-regular if $t=s$;
  \item [{\rm (3)}]\ $\Sigma_{3,t,s}$ is $2$-arc-regular if $t=s+1$, and $(t, s)\neq(2, 1)$;
  \item [{\rm (4)}]\ $\Sigma_{p,t,s}$ is $1$-arc-regular if $p^{t-s}>3$.
\end{enumerate}
\end{theorem}

\f\demo By Magma~\cite{Magma}, we can obtain $(1)$. If $(p,t,s)=(3,1,1)$ then by Magma~\cite{Magma},
we can show that $\Sigma_{3,1,1}$ is $2$-arc-regular. In what follows, we assume that $(p,t,s)\neq (3,2,1), (3,1,1)$.

Set $\G=\Sigma_{p,t,s}$ and $H=\H_{p,t,s}$. We shall first prove that $\G$ is a normal bi-Cayley graph over $H$.
By Lemma~\ref{57}, we may assume that $p=3$.
Since $(p,t,s)\neq (3,1,1), (3,2,1)$, one has $|H|=3^{t+s+1}\geq 3^5$. Let $n=t+s+1$.

Let $A=\Aut(\G)$ and let $P$ be a Sylow $3$-subgroup of $A$ such that $R(H)\leq P$. Then $R(H)\unlhd P$. 
By Lemma~\ref{Example2}, $\G$ is symmetric. We first prove the following claim.\medskip

\f{\bf Claim 1}\ $P\unlhd A$.\medskip

Let $M\unlhd A$ be maximal subject to that $M$ is intransitive on both $H_0$ and $H_1$. By Proposition~\ref{3orbits}, $M$ is semiregular on $V(\G)$ and the quotient graph $\G_M$ of $\G$ relative to $M$ is a cubic graph with $A/M$ as an arc-transitive group of automorphisms. Assume that $|M|=3^\ell$. Then $|V(\G_M)|=2\cdot 3^{n-\ell}$. If $n-\ell\leq 3$, then by \cite{Sym768}, $\G_M$ is isomorphic to $\F006A$, $\F018A$ or $\F054A$, and then by Magma~\cite{Magma}, $\Aut(\G_M)$ has a normal Sylow $3$-subgroup. It follows that $P/M\unlhd A/M$, and so $P\unlhd A$, as claimed.

Now suppose that $n-\ell>3$. Take a minimal normal subgroup $N/M$ of $A/M$. By Proposition~\ref{N}, $N/M$ is an elementary abelian $3$-group. By the maximality of $M$, $N$ is transitive on at least one of $H_0$ and $H_1$, and so $3^{n}\ |\ |N|$. If $3^{n+1}\ |\ |N|$, then $P=N\unlhd A$, as claimed. Now assume that $|N|=3^n$. We have $N$ is transitive on both $H_0$ and $H_1$. Then $N$ is semiregular on both $H_0$ and $H_1$, and then $\G_M$ would be a cubic bi-Cayley graph on $N/M$. Since $\G_M$ is connected, by Proposition~\ref{bicayley}, $N/M$ is generated by two elements, and so $N/M\cong\mz_3$ or $\mz_3\times\mz_3$. This implies that $|V(\G_M)|=6$ or $18$, contrary to the assumption that $|V(\G_M)|=2\cdot 3^{n-\ell}>18$, completing the proof of our claim. \medskip

By Claim~1, we have $P\unlhd A$. Since $|P: R(H)|=3$, one has $\Phi(P)\leq R(H)$. As $H$ is non-abelian, one has $\Phi (P)<R(H)$ for otherwise, we would have $P$ is cyclic and so $H$ is cyclic which is impossible. Then $\Phi(P)$ is intransitive on both $H_0$ and $H_1$, the two orbits of $R(H)$ on $V(\G)$. Since $\Phi(P)$ is characteristic in $P$, $P\unlhd A$ gives that $\Phi(P)\unlhd A$. By Propositions~\ref{3orbits}, the quotient graph $\G_{\Phi(P)}$ of $\G$ relative to $\Phi(P)$ is a cubic graph with $A/\Phi(P)$ as an arc-transitive group of automorphisms. Furthermore, $P/\Phi(P)$ is transitive on the edges of $\G_{\Phi(P)}$. Since $P/\Phi(P)$ is abelian, it is easy to see that $\G_{\Phi(P)}\cong K_{3,3}$, and so $P/\Phi(P)\cong\mz_3\times\mz_3$.

Let $\Phi_2$ be the Frattini subgroup of $\Phi(P)$.
Then $\Phi_2\unlhd A$ because $\Phi_2$ is characteristic in $\Phi (P)$ and $\Phi(P)\unlhd A$.
Let $\Phi_3$ be the Frattini subgroup of $\Phi_2$.
Similarly, we have $\Phi_3\unlhd A$.
Now we prove the following claim. \medskip

\f{\bf Claim 2}\ $\Phi(P)/\Phi_2\cong\z_3\times\z_3\times\z_3$ and $\Phi_2/\Phi_3\cong\z_3\times\z_3$.\medskip

Since $P/\Phi(P)\cong\mz_3\times\mz_3$ and $|P:R(H)|=3$, we have $|R(H):\Phi(P)|=3$.
Then $\Phi(P)$ is a maximal subgroup of $R(H)$. And then by Lemma~\ref{comput2}~(5), we have $\Phi(P)$ is isomorphic to one of the following four groups:
\begin{align*}
 & M_1=\lg a\rg\times\lg b^3\rg\times\lg c\rg, ~~~~~~~~~~~~~M_2=\lg a^3\rg\times\lg b\rg\times\lg c\rg, \\
 & M_3=\lg ab\rg\times\lg b^3\rg\times\lg c\rg, ~~~~~~~~~~~M_4=\lg ab^{-1}\rg\times\lg b^3\rg\times\lg c\rg.
\end{align*}
It follows that $\Phi(P)/\Phi_2\cong\z_3\times\z_3\times\z_3$. Then $\Phi_2$ is isomorphic to one of the following four groups:
$$Q_1=\lg a^3\rg\times\lg b^9\rg\cong\z_{3^{t-1}}\times\z_{3^{s-2}},~~~~~~~~~ Q_2=\lg a^9\rg\times\lg b^3\rg\cong\z_{3^{t-2}}\times\z_{3^{s-1}},~~~$$
$$Q_3=\lg a^3b^3\rg\times\lg b^9\rg\cong\z_{3^{t-1}}\times\z_{3^{s-2}},~~~~~ Q_4=\lg a^{3}b^{-3}\rg\times\lg b^9\rg\cong\z_{3^{t-1}}\times\z_{3^{s-2}}.$$

It implies that $\Phi_2/\Phi_3\cong\z_3\times\z_3$, as claimed.\medskip

Clearly, $\Phi_3\leq \Phi(P)<R(H)$, so $\Phi_3$ is intransitive on both $H_0$ and $H_1$.
Consider the quotient graph $\G_{\Phi_3}$ of $\G$ relative to $\Phi_3$. By Propositions~\ref{3orbits}, $\G_{\Phi_3}$ is a cubic graph with $A/\Phi_3$ as an arc-transitive group of automorphisms. Furthermore, $\G_{\Phi_3}$ is a bi-Cayley graph over the group $R(H)/\Phi_3$ of order $2\cdot 3^6$.
Then by \cite{Conder}, $\G_{\Phi_3}\cong \C1458.1$, $\C1458.2$, $\C1458.3$, $\C1458.4$, $\C1458.5$, $\C1458.6$, $\C1458.7$, $\C1458.8$, $\C1458.9$, $\C1458.10$ or $\C1458.11$. By Magma \cite{Magma}, if $\G_{\Phi_3}\cong \C1458.1$, $\C1458.3$, $\C1458.4$, $\C1458.8$, $\C1458.9$, $\C1458.10$ or $\C1458.11$, then $\Aut(\G_{\Phi_3})$ does not have an abelian or inner-abelian semiregular subgroup of order $729$, a contradiction.
If $\G_{\Phi_3}\cong \C1458.2$, $\C1458.5$, $\C1458.6$ or $\C1458.7$, then by Magma \cite{Magma}, all semiregular subgroups of $\Aut(\G_{\Phi_3})$ of order $729$ are normal, and so $R(H)/\Phi_3\unlhd \Aut(\G_{\Phi_3})$. It follows that $R(H)/\Phi_3\unlhd A/\Phi_3$, and so $R(H)\unlhd A$.

By now we have shown that $\Sigma_{p,t,s}$ is normal. By \cite[Theorem~1.1]{NET}, $\Sigma_{p,t,s}$ is at most $2$-arc-transitive. Recall that Lemma~\ref{Example2} already proved that $\Sigma_{p,t,s}$ is at least $1$-arc-transitive.

Let $t=s$. In this case, we have $k=0$ and \[\Sigma_{p, t, t}=\BiCay(\H_{p, t, t},\emptyset, \emptyset, \{1, a, b\}.\]
It is easy to see that $\H_{p,t,t}$ has an automorphism $\g$ swapping $a$ and $b$. Then $\s_{\g,1}\in\Aut(\Sigma_{p,t,t})_{1_{0}1_1}$ and $\s_{\g,1}$
swaps $a_1$ and $b_1$. Thus, $\Sigma_{p,t,t}$ is $2$-arc-regular.

Let $t=s+1$ and $p=3$. In this case, we have $k^2-k+1\equiv0\ (\mod 3)$ and so $k=2$ since $k\in\mz_3^*$.
Then \[\Sigma_{p, s+1, s}=\BiCay(\H_{p, s+1, s},\emptyset, \emptyset, \{1, a, ba^2\}.\]
By Lemma~\ref{comput2}~(1), we see that $(ba^2)^2=b^2a^4c^2$, and so $a(ba^2)^{-2}=a^{-3}b^{-2}c^{-2}$, which has the same order as $b$.
Noticing that $o(a)=o(ba^2)$, by Lemma~\ref{comput2}~(4), $\H_{p,s+1,s}$ has an automorphism $\g$ taking $(a, b)$ to $(ba^2, a(ba^2)^{-2})$.
Furthermore, $\g$ swaps $a$ and $ba^2$. Then $\s_{\g,1}\in\Aut(\Sigma_{p,s+1,s})_{1_{0}1_1}$ and $\s_{\g,1}$
swaps $a_1$ and $(ba^2)_1$. Thus, $\Sigma_{p,s+1,s}$ is $2$-arc-regular.

Let $p^{t-s}>3$. In this case, $\Sigma_{p, t, s}=\BiCay(\H_{p, t, s},\emptyset, \emptyset, \{1, a, ba^{k}\}),$
where $k\in\mz_{p^{t-s}}^*$ satisfies $k^2-k+1\equiv0\ (\mod p^{t-s})$. If $\Sigma_{p,t,s}$ is $2$-arc-regular, then by \cite[Theorem~1.1]{NET}, $\H_{p,t,s}$ has an automorphism $\g$ swapping $a$ and $ba^k$, and then $b^\g=(ba^k)^\g(a^{-k})^\g=a(ba^k)^{-k}$. It follows that $1=(a(ba^k)^{-k})^{p^s}=(a^{p^s})^{1-k^2}$, and so $1-k^2\equiv 0\ (\mod p^{t-s})$. Combining this with the equation $k^2-k+1\equiv0\ (\mod p^{t-s})$, we have $k\equiv 2\ (\mod p^{t-s})$, forcing $p^{t-s}=3$, a contradiction. Thus, $\Sigma_{p, t, s}$ is $1$-arc-regular.\hfill\qed

\section{Proof of Theorem~\ref{classify-all}}

The goal of this section is to prove Theorem~\ref{classify-all}.\medskip


\f{\bf Proof of Theorem~\ref{classify-all}}\ To complete the proof, by Theorem~\ref{s-regular},
it suffices to prove that every cubic edge-transitive bi-Cayley graph over $\H_{p,t,s}$ is isomorphic to $\Sigma_{p,t,s}$.

Let $H=\H_{p,t,s}$, and let $\G=\BiCay(H, \R, \L, S)$ be a connected cubic edge-transitive bi-Cayley graph over $\H_{p,t,s}$
Set $A=\Aut(\G)$. By Lemma~\ref{57}, we have $\G$ is normal edge-transitive. It follows that the two orbits $H_0, H_1$ of $R(H)$ on $V(\G)$ do not contain edges of $\G$, and so $\R=\L=\emptyset$. By Proposition~\ref{bicayley}, we may assume that $S=\{1, x, y\}$ for $x, y\in H$. Since $\G$ is connected, by Proposition~\ref{bicayley}, we have $H=\lg S\rg=\lg x, y\rg$.

Since $\G$ is normal edge-transitive, by Proposition~\ref{bicayleyaut}, there exists $\sigma_{\a, h}\in A_{1_0}$, where $\a\in\Aut(H)$ and $h\in H$, such that $\s_{\a, h}$ cyclically permutates the three elements in $\G(1_0)=\{1_1, x_1, y_1\}$. Without loss of generality, assume that $(\sigma_{\a, h})_{|\G(1_0)}=(1_1\ x_1\ y_1)$. Then $x_1=(1_1)^{\sigma_{\a, h}}=h_1$, implying that $x=h$.
Furthermore, $y_1=(x_1)^{\sigma_{\a, h}}=(xx^{\a})_1$ and $1_1=(y_1)^{\sigma_{\a, h}}=(xy^{\a})_1$.
It follows that $x^{\a}=x^{-1}y$ and $y^{\a}=x^{-1}$.

Recall that
$$H=\H_{p,t,s}=\lg a, b, c\ |\ a^{p^t}=b^{p^s}=c^p=1, [a, b]=c, [c, a]=[c, b]=1\rg,$$
where $t\geq s\geq 1$. We first prove the following claim.\medskip

\f{\bf Claim}\ $o(x)=o(y)=o(x^{-1}y)=p^{t}$ and $\lg x^{p^s}\rg=\lg y^{p^s}\rg$. \medskip

Since $x^{\a}=x^{-1}y$ and $y^{\a}=x^{-1}$, we have $o(x)=o(y)=o(x^{-1}y)$. Denote by ${\rm exp}(H)$ the exponent of $H$.
Since $H=\lg x, y\rg$, by Lemma~\ref{comput2}~(3), we have $o(x)=o(y)=o(x^{-1}y)={\rm exp}(H)=p^t$.
Note that $H=\lg a, b\rg=\lg x, y\rg$. Again by Lemma~\ref{comput2}~(3), we have $\lg x^{p^s}\rg\leq \lg a\rg$ and $\lg y^{p^s}\rg\leq \lg a\rg$.
Since $o(x)=o(y)=p^{t}$, we have $\lg x^{p^s}\rg=\lg y^{p^s}\rg$, as claimed.\medskip

Now we are ready to finish the proof. If $t=s$, then by Lemma~\ref{comput2}~(4), there exists an automorphism of $H$
sending $(x, y)$ to $(a,b)$, and by Proposition~\ref{bicayley}~(3), we have $\G\cong\Sigma_{p,t,t}$.

Suppose now that $t>s$. By Claim, we have $\lg x^{p^s}\rg=\lg y^{p^s}\rg$.
Then there exists $k\in\z_{p^{t-s}}^*$ such that $y^{p^s}=x^{kp^s}$, and so $(yx^{-k})^{p^s}=1$.
So $o(yx^{-k})=o(b)=p^s$. By Claim, we have $o(x)=o(a)=p^t$. Since $H=\lg x, y\rg=\lg x, yx^{-k}\rg$,
by Lemma~\ref{comput2}~(4), there exists $\gamma\in\Aut(H)$ such that $a^{\gamma}=x$ and $b^{\gamma}=yx^{-k}$. It follows that
$$H=\lg x, yx^{-k}, z\ |\ x^{p^t}=(yx^{-k})^{p^s}=z^p=1, [x, yx^{-k}]=z, [z, x]=[z, yx^{-k}]=1\rg,$$
and $S=\{1, x, y\}=\{1, x, (yx^{-k})x^k\}$. Clearly, $S^{\g^{-1}}=\{1, a, ba^k\}$. By Proposition~\ref{bicayley}~(3),
we may assume that $\G=\BiCay(H, \emptyset, \emptyset, \{1, a, ba^k\})$.

Since $\G$ is normal edge-transitive, by Proposition~\ref{bicayleyaut}, there exists $\sigma_{\theta, g}\in \Aut(\G)_{1_0}$,
where $\theta\in\Aut(H)$ and $g\in H$, such that $\s_{\theta, g}$ cyclically permutates the three elements in $\G(1_0)=\{1_1, a_1, (ba^k)_1\}$.
Without loss of generality, assume that $(\sigma_{\theta, g})_{|\G(1_0)}=(1_1\ a_1\ (ba^k)_1)$. Then $a_1=(1_1)^{\sigma_{\theta, g}}=g_1$, implying that $a=g$.
Furthermore, we have $$(ba^k)_1=(a_1)^{\sigma_{\theta, g}}=(aa^{\theta})_1, ~~1_1=(ba^k)_1^{\sigma_{\theta, g}}=(a(ba^k)^{\theta})_1. $$
Then
$$a^{\theta}=a^{-1}ba^k=ba^{k-1}c^{-1}, ~~b^{\theta}=a^{-1}{(a^{\theta}})^{-k}=a^{-1}(ba^{k-1}c^{-1})^{-k}. $$
This implies that $o(a^{-1}(ba^{k-1}c^{-1})^{-k})=o(b)=p^s$.
By Lemma~\ref{comput2}~(1) and (3), we have $o(a^{k-1})=p^t$ and $(a^{-1}(ba^{k-1}c^{-1})^{-k})^{p^s}=a^{-(k^2-k+1)p^s}=1$.
It follows that $k^2-k+1\equiv0\ (p^{t-s})$, and hence $\G\cong\Sigma_{p,t,s}$. \hfill\qed

\medskip
\f {\bf Acknowledgements:}\ This work was partially supported by the National
Natural Science Foundation of China (11271012) and the Fundamental
Research Funds for the Central Universities (2015JBM110).

\end{document}